\documentclass[11pt]{amsart}

\usepackage[margin=1.25in]{geometry}

\usepackage[dvips]{graphicx}
\usepackage{calc}
\usepackage{color}
\usepackage{amsmath}
\usepackage{amssymb}
\usepackage{amscd}
\usepackage{amsthm}
\usepackage{amsbsy}
\usepackage{delarray}
\usepackage{enumerate}
\usepackage[T1]{fontenc}
\usepackage{inputenc}
\usepackage{enumerate}
\usepackage{hyperref}

\begin{document}

\setlength{\parindent}{5mm}
\renewcommand{\leq}{\leqslant}
\renewcommand{\geq}{\geqslant}
\newcommand{\N}{\mathbb{N}}

\newcommand{\Z}{\mathbb{Z}}
\newcommand{\R}{\mathbb{R}}
\newcommand{\C}{\mathbb{C}}
\newcommand{\F}{\mathbb{F}}
\newcommand{\g}{\mathfrak{g}}
\newcommand{\h}{\mathfrak{h}}
\newcommand{\K}{\mathbb{K}}
\newcommand{\RN}{\mathbb{R}^{2n}}
\newcommand{\ci}{c^{\infty}}
\newcommand{\derive}[2]{\frac{\partial{#1}}{\partial{#2}}}
\renewcommand{\S}{\mathbb{S}}
\renewcommand{\H}{\mathbb{H}}
\newcommand{\eps}{\varepsilon}
\renewcommand{\r} {\textbf{r}}
\renewcommand{\L}{\mathcal{L}}

\theoremstyle{plain}
\newtheorem{theo}{Theorem}
\newtheorem{prop}[theo]{Proposition}
\newtheorem{lemma}[theo]{Lemma}
\newtheorem{definition}[theo]{Definition}
\newtheorem*{notation*}{Notation}
\newtheorem*{notations*}{Notations}
\newtheorem{corol}[theo]{Corollary}
\newtheorem{conj}[theo]{Conjecture}
\newtheorem{question}[theo]{Question}
\newtheorem*{question*}{Question}

\newenvironment{demo}[1][]{\addvspace{8mm} \emph{Proof #1.
    ---~~}}{~~~$\Box$\bigskip}

\newlength{\espaceavantspecialthm}
\newlength{\espaceapresspecialthm}
\setlength{\espaceavantspecialthm}{\topsep} \setlength{\espaceapresspecialthm}{\topsep}

\newenvironment{example}[1][]{\refstepcounter{theo} 
\vskip \espaceavantspecialthm \noindent \textsc{Example~\thetheo
#1.} }%
{\vskip \espaceapresspecialthm}

\newenvironment{remark}[1][]{\refstepcounter{theo} 
\vskip \espaceavantspecialthm \noindent \textsc{Remark~\thetheo
#1.} }%
{\vskip \espaceapresspecialthm}

\def\Homeo{\mathrm{Homeo}}
\def\Hameo{\mathrm{Hameo}}
\def\Diffeo{\mathrm{Diffeo}}
\def\Symp{\mathrm{Symp}}
\def\Id{\mathrm{Id}}
\newcommand{\norm}[1]{||#1||}
\def\Ham{\mathrm{Ham}}
\def\Hamtilde{\widetilde{\mathrm{Ham}}}
\def\Crit{\mathrm{Crit}}
\def\Spec{\mathrm{Spec}}
\def\osc{\mathrm{osc}}
\def\Cal{\mathrm{Cal}}

\title[]{The displaced disks problem via symplectic topology} 
\author{Sobhan Seyfaddini}
\date{\today}

\address{SS: D\'epartement de Math\'ematiques et Applications de l'\'Ecole Normale Sup\'erieure, 45 rue d'Ulm, F 75230 Paris cedex 05}
\email{sobhan.seyfaddini@ens.fr}


\begin{abstract} 
We prove that a $C^0$--small area preserving homeomorphism of a closed surface with vanishing mass flow can not displace a topological disk of large area.  This resolves the displaced disks problem posed by F. B\'eguin, S. Crovisier, and F. Le Roux. 
\end{abstract}

\maketitle

\section{Introduction}
  Let $\Sigma$ be a closed surface endowed with an area form $\Omega$.  We denote by $Homeo^{\Omega}(\Sigma)$ the group of area preserving homeomorphisms of $\Sigma$, by $Homeo_0^{\Omega}(\Sigma)$ the path component of identity in $Homeo^{\Omega}(\Sigma)$, and by $Ham(\Sigma)$ the group of Hamiltonian diffeomorphisms of $\Sigma$.  We will be studying, $\mathcal{G},$ the $C^0$ closure of $Ham(\Sigma)$ inside $Homeo_0^{\Omega}(\Sigma)$. 

	The group $\mathcal{G}$ is a well-studied dynamical object: it is precisely the set of elements of $Homeo_0^{\Omega}(\Sigma)$ with vanishing \emph{mass flow}.  For a definition of the mass flow homomorphism, which is also known as the \emph{mean rotation vector}, see Section $5$ of \cite{Fa}.  Equivalently,  $\mathcal{G}$ can be described as the set of elements of $Homeo_0^{\Omega}(\Sigma)$ with zero flux; see Appendix A.5 of \cite{Fa}.   In this note, we only work with the description of $\mathcal{G}$ as the $C^0$ closure of $Ham(\Sigma)$.  It is well-known that in the case of $S^2$, $\mathcal{G} = Homeo_0^{\Omega}(S^2)$.
	
	\medskip
  Recall that a homeomorphism $\phi$ is said to displace a set $B$ if $\phi(B) \cap B =\emptyset$.   For $a > 0$, define $\mathcal{G}_a = \{ \theta \in \mathcal{G}:\theta \text{ displaces a topological disk of area at least } a \}.$   The displaced disks problem, posed by F. B\'eguin, S. Crovisier, and F. Le Roux, asks the following.
\begin{question*}(B\'eguin, Crovisier, Le Roux)
Does the $C^0$ closure of $\mathcal{G}_a$ contain the identity?
\end{question*} 

The initial motivation of  B\'eguin, Crovisier, and Le Roux for posing this beautiful question is as follows: $\mathcal{G}$ is a normal subgroup of $Homeo^{\Omega}(\Sigma).$ B\'eguin, Crovisier, and Le Roux were interested in knowing whether the conjugacy class of an element of $\mathcal{G}$ could be $C^0$--dense in all of $\mathcal{G}$.    As we will see below, a negative answer to the displaced disks problem provides a negative answer to this question as well.  

\medskip

  Equip $\Sigma$ with a Riemannian distance $d$ and define the $C^0$ distance between two homeomorphisms $\phi, \, \psi$ by $d_{C^0}(\phi, \psi) :=  \max\{ \max_{x \in \Sigma} d(\phi(x), \psi(x)), \max_{x \in \Sigma} d(\phi^{-1}(x), \psi^{-1}(x))\}$.  Our main result provides a negative answer to the displaced disks problem.

\begin{theo} \label{main_theo}
 There exists $\epsilon > 0$ such that $d_{C^0}(Id, \theta ) \geq \epsilon$ for all $\theta \in \mathcal{G}_a$.
\end{theo}
 
Let $C(\phi) = \{ \psi \phi \psi^{-1} : \psi \in Homeo^{\Omega}(\Sigma) \}$ denote the conjugacy class of $\phi \in \mathcal{G}$.  Observe that $\mathcal{G}_a$ is invariant under conjugation by area preserving homeomorphisms and hence, if $\phi \neq Id$, then  $C(\phi) \subset \mathcal{G}_a$, for some $a>0$.  Theorem \ref{main_theo} immediately yields the next result.
\begin{corol}\label{cor}
 $C(\phi)$ is not $C^0$--dense in $\mathcal{G}$.  
\end{corol}

The following facts were pointed out to me by B\'eguin,  Crovisier, and Le Roux.

\noindent \textbf{1. }Theorem \ref{main_theo}  would not hold without requiring  $\theta \in \mathcal{G}$.  Indeed, it is not difficult to see that a $C^0$--small translation of the torus displaces a disk with area nearly equal to half the total area of the torus.

\noindent \textbf{2. }Corollary \ref{cor} does not hold for arbitrary elements of $Homeo_0(\Sigma)$.  See Remark 7.11 of \cite{Gu} for an example of a homeomorphism of $S^2$ whose conjugacy class is $C^0$--dense in $Homeo_0(S^2)$.  Once it is established that the conjugacy class of a homeomorphism is dense, then it is easy to see  that the conjugacy class of that homeomorphism is a $G_\delta$ set.  Hence, we conclude that the conjugacy class of a generic homeomorphism of $S^2$ is dense.

\noindent \textbf{3. }Suppose that $\Sigma \neq S^2$.  It follows from the work of Gaumbado and Ghys \cite{GG} and Entov, Polterovich, and Py \cite{EPP} that $\mathcal{G}$ carries $C^0$--continuous and homogeneous quasimorphisms; see Theorem 1.2 of \cite{EPP}.  Corollary \ref{cor} follows immediately as homogeneous quasimorphisms are constant on conjugacy classes.  

\section{Proof of Theorem \ref{main_theo}}
Our proof uses Floer theoretic invariants of Hamiltonian diffeomorphisms.  In particular, we use the theory of spectral invariants, or action selectors, introduced by C. Viterbo, M. Schwarz, and Y.-G. Oh \cite{viterbo1, schwarz, oh2}.  An important consequence of this theory is that the group of Hamiltonian diffeomorphisms of a closed symplectic manifold $M$ admits a conjugation invariant norm $\gamma : Ham(M) \rightarrow [0, \infty)$ \footnote{$\gamma$ is usually defined on the universal cover of $Ham$.  For $\phi \in Ham$, one can define $\gamma(\phi)$ by taking infimum over all paths which end at $\phi$.}.  Being a conjugation invariant norm, $\gamma$ satisfies the following axioms:

\begin{enumerate}[i.]
\item \label{positivity} $\gamma(\phi) \geq 0$ with equality if and only if $\phi = Id$,
\item \label{inverse} $\gamma(\phi) = \gamma(\phi^{-1})$,
\item \label{triangle} $\gamma(\phi \psi) \leq \gamma(\phi) + \gamma(\psi)$,
\item \label{invariance}$\gamma(\psi \phi \psi^{-1})  = \gamma (\phi)$.
\end{enumerate}

An important feature of $\gamma$ is the fact that it satisfies the so-called \emph{energy-capacity} inequality.  In the case of a closed surface $\Sigma$ the energy-capacity inequality states that if $\phi \in Ham(\Sigma)$ displaces a disk of area $a$ then 
\begin{equation} \label{eq: energy-capacity}
a \leq \gamma(\phi).
\end{equation}
Theorem 2 of \cite{seyfaddini11} provides the final step of our solution.  According to this theorem , for a closed surface $\Sigma$, of genus $g$, there exist constants $C, \delta >0$ such that $\forall \phi \in Ham(\Sigma)$ if $d_{C^0}(Id, \phi) \leq \delta$ then 
\begin{equation}\label{Lip_Continuity}
\gamma(\phi) \leq C \, d_{C^0}(Id, \phi)^{2^{-2g-1}}.
\end{equation}

\medskip
We now prove Theorem \ref{main_theo}.  For a contradiction, suppose it does not hold and pick a sequence $\theta_i \in \mathcal{G}_a$ which converges uniformly to the identity and conclude from Inequality (\ref{Lip_Continuity}) that $\gamma(\theta_i) \to 0$. But this is impossible because the energy-capacity inequality (\ref{eq: energy-capacity}) implies that $\gamma|_{\mathcal{G}_a}  \geq a.$

\subsection{Extending $\gamma$ to $\mathcal{G}$}  We will finish this note by showing that the conjugation invariant norm $\gamma$ extends continuously to $\mathcal{G}.$ We need a small modification of Inequality (\ref{Lip_Continuity}).  Let $C, \delta$ denote the constants appearing in this inequality and  suppose that $d_{C^0}(\psi, \phi) \leq \delta$, where  $\psi, \phi \in Ham(\Sigma)$.  Applying Inequality (\ref{Lip_Continuity}) to $\psi \phi^{-1}$ we obtain that $\gamma(\psi \phi^{-1}) \leq C \, d_{C^0}(Id, \psi \phi^{-1})^{2^{-2g-1}} \leq C \, d_{C^0}(\phi, \psi)^{2^{-2g-1}}$; the latter inequality follows from the definition of $d_{C^0}$.  Now, using Axiom (\ref{triangle}) of $\gamma$ we see that $\gamma(\psi) - \gamma(\phi) \leq \gamma(\psi \phi^{-1})$.   Hence, $\gamma(\psi) - \gamma(\phi) \leq C \, d_{C^0}(\psi, \phi)^{2^{-2g-1}}.$  Similarly, we obtain the same upper bound for $\gamma(\phi) - \gamma(\psi)$.  Therefore, we have proven that  $\forall \psi, \, \phi \in Ham(\Sigma)$ if $d_{C^0}(\psi, \phi) \leq \delta$ then 
$$|\gamma(\psi) - \gamma(\phi)| \leq C \, d_{C^0}(\psi, \phi)^{2^{-2g-1}}.$$

We see that $\gamma$ is uniformly continuous with respect to $d_{C^0}$ and so, it extends continuously to $\mathcal{G}$.  Clearly, the extension $\gamma: \mathcal{G} \rightarrow \mathbb{R}$ satisfies Inequalities (\ref{eq: energy-capacity}) and (\ref{Lip_Continuity}), in addition to the four stated axioms of conjugation invariant norms.  

\subsection*{Aknowledgements}
 The problem solved by Corollary \ref{cor} was brought to my attention by Dmitri Burago, Sergei Ivanov, and Leonid Polterovich .  Furthermore,  Leonid Polterovich suggested that my results from \cite{seyfaddini11} would lead to a solution.   I am very grateful to the three of them. 
I would like to thank Fr\'ed\'eric Le Roux for organizing a set of lectures on the contents of this note.  While preparing this article I learned a great deal about the circle of ideas surrounding the displaced disks problem from  Fran\c{c}ois B\'eguin, Sylvain Crovisier, and Fr\'ed\'eric Le Roux; I am thankful to them for their invaluable input.  Lastly, I would like to thank Vincent Humili\`ere, R\'emi Leclercq, Leonid Polterovich and Alan Weinstein for their helpful comments on preliminary drafts of this note.

\bibliographystyle{abbrv}
\bibliography{biblio}

\begin{thebibliography}{1}

\bibitem{EPP}
M.~Entov, L.~Polterovich, and P.~Py.
\newblock On continuity of quasimorphisms for symplectic maps.
\newblock In {\em Perspectives in analysis, geometry, and topology}, volume 296
  of {\em Progr. Math.}, pages 169--197. Birkh\"auser/Springer, New York, 2012.
\newblock With an appendix by Michael Khanevsky.

\bibitem{Fa}
A.~Fathi.
\newblock Structure of the group of homeomorphisms preserving a good measure on
  a compact manifold.
\newblock {\em Ann. Sci. \'Ecole Norm. Sup. (4)}, 13(1):45--93, 1980.

\bibitem{GG}
J.-M. Gambaudo and {\'E}.~Ghys.
\newblock Commutators and diffeomorphisms of surfaces.
\newblock {\em Ergodic Theory Dynam. Systems}, 24(5):1591--1617, 2004.

\bibitem{Gu}
P.-A. Guih{\'e}neuf.
\newblock {\em Propri\'et\'es dynamiques g\'en\'eriques des hom\'eomorphismes
  conservatifs}, volume~22 of {\em Ensaios Matem\'aticos [Mathematical
  Surveys]}.
\newblock Sociedade Brasileira de Matem\'atica, Rio de Janeiro, 2012.

\bibitem{oh2}
Y.-G. {Oh}.
\newblock Construction of spectral invariants of {H}amiltonian paths on closed
  symplectic manifolds.
\newblock In {\em The breadth of symplectic and {P}oisson geometry}, volume 232
  of {\em Progr. Math.}, pages 525--570. Birkh\"auser Boston, Boston, MA, 2005.

\bibitem{schwarz}
M.~Schwarz.
\newblock On the action spectrum for closed symplectically aspherical
  manifolds.
\newblock {\em Pacific J. Math.}, 193(2):419--461, 2000.

\bibitem{seyfaddini11}
S.~{Seyfaddini}.
\newblock ${C}^0$--limits of {H}amiltonian paths and the {O}h--{S}chwarz
  spectral invariants.
\newblock {\em Int. Math. Res. Not.}
\newblock (to appear).

\bibitem{viterbo1}
C.~Viterbo.
\newblock Symplectic topology as the geometry of generating functions.
\newblock {\em Math. Annalen}, 292:685--710, 1992.

\end{thebibliography}

\end{document}